# Subnormal subgroups of groups and ideals of group rings.

Scherback Denis

***Abstract***. It is known that the group of Hamiltonian quaternions $\mathbb{H}$ is not abelian. This article aims to show that all subgroups of $\mathbb{H}$ are abelian and group rings $\mathbb{Z}_2[\mathbb{H}]$ and $\mathbb{Z}_3[\mathbb{H}]$ are not fields. And also finding all the ideals of group rings $\mathbb{Z}_2[\mathbb{H}]$ and $\mathbb{Z}_3[\mathbb{H}]$.

***Definition***: If the ring $K$ is finite-dimensional vector space over a field $P$ and satisfying the condition of associativity

$$(\alpha u)v = u(\alpha v) = \alpha(uv), \quad \alpha \in P$$

then it is called associative algebra over the field $P$.
If we exclude the condition of associativity we get the general concept of (linear) algebra.
Suppose $K$ a 4-dimensional vector space with basis elements $e, i, j, k$
We define equality

$$i^2 = -e\alpha, \quad j = -e\beta$$

where $\alpha, \beta$ is an arbitrary elements from $P$.

The resulting algebra $K$ is called associative algebra of quaternions, its components look like

$$x = ex_0 + ix_1 + jx_2 + kx_3$$

If $\alpha = \beta = 1$, then we obtain the Hamiltonian quaternions.
Check associativity of the algebra of quaternions.
To do this, check the 27 equations of the form $a(bc) = (ab)c$. We can avoid this by setting the isomorphism of the algebra of quaternions over the field $\mathbb{R}$ and an algebra of matrices of special form over $\mathbb{C}$.
Namely:

$$e \leftrightarrow E = \begin{pmatrix} 1 & 0 \\ 0 & 1 \end{pmatrix}, \quad i \leftrightarrow I = \begin{pmatrix} i & 0 \\ 0 & -i \end{pmatrix}$$

$$j \leftrightarrow J = \begin{pmatrix} 0 & 1 \\ -1 & 0 \end{pmatrix}, \quad k \leftrightarrow K = \begin{pmatrix} 0 & i \\ i & 0 \end{pmatrix}$$

Equality $I^2 = J^2 = K^2 = -E$, $IJ = -JI = K$, $JK = -KJ = I$, $KI = -IK = J$, it is easy to verify. This means that the space of matrices formed by the matrices $E, I, J, K$ forms an algebra which is isomorphic to the algebra of quaternions. From the associativity of matrix multiplication we conclude about the associativity of the algebra of quaternions.

## Properties of group of Hamiltonian quaternions $\mathbb{H}$

Consider the group of Hamiltonian quaternions
$$\mathbb{H} = \{e, -e, i, -i, j, -j, k, -k\}$$
it is an Abelian group by virtue of a prescribed multiplication of the elements $i, j, k$.
Consider the subgroup $\mathbb{H}$

$\{e\}$
$\{e, -e\}$
$\{e, e, i, -i\}$
$\{e, e, j, -j\}$
$\{e, e, k, -k\}$

Denote
$A_e = \{e, -e\}$
$A_i = \{e, -e, i, -i\}$
$A_j = \{e, -e, j, -j\}$
$A_k = \{e, -e, k, -k\}$

***Statement 1***: All subgroups of $\mathbb{H}$ are abelian.
Really
$\{e\} = <e>$
$\{e, -e\} = <-e>$
$\{e, e, i, -i\} = <-i>$
$\{e, -e, j, -j\} = <-j>$
$\{e, -e, k, -k\} = <-k>$
All subgroups of $\mathbb{H}$ are cyclic, and hence abelian.
Also $A_e = \{e, -e\}$ is the center of the group $\mathbb{H}$.

***Statement 2***: The group of $\mathbb{H}$ is nilpotent.
***Definition***: Suppose that the group $G$ has an invariant series
$$E = A_0 \subset A_1 \subset A_2 \subset ... \subset A_n = G$$
This series is called a central series, if $t = 0, 1, ..., n-1$ the factor group $A_{t+1}/A_t$ is in the center factor groups, in other words, if the mutual commutator
$$[A_{t+1} \quad G] \subset A_t, \quad t = 0, 1, ..., n-1$$
We prove that the series $\mathbb{H} \geq A_i \geq A_e \geq \{e\}$ is central.

Consider $[A_e \mathbb{H}]$
Let
$a \in A_e$
$h \in \mathbb{H}$
$a^{-1}h^{-1}ah = a^{-1}ah^{-1}h = e \quad \forall \, a \in A_e, \, h \in \mathbb{H}$
so $[A_e \mathbb{H}] \in \{e\}$.

Next, we consider $[A_i \mathbb{H}]$.
1) If $a \in A_e$ or $h \in A_e$ then
$$[ah] = e$$
2) $a \notin A_e$ и $h \notin A_e$
$$a^{-1}h^{-1}ah = a^{-1}ah^{-1}h = -e$$
so $[A_i \mathbb{H}] \subset \{e\}$.

Analogously for $[\mathbb{H}\mathbb{H}]$ we can prove that $[\mathbb{H}\mathbb{H}] \subset A_e \subset A_i$.
Thus we have proved that the series $\mathbb{H} \geq A_i \geq A_e \geq \{e\}$ is central.

So $\mathbb{H}$ is nilpotent group.

***Definition***: The finite normal or invariant series of the group is called solvable series if all its factors are abelian.

The group $G$ is called solvable if it satisfies any of the next conndition:
1) The group $G$ has a finite solvable normal series.
2) The group $G$ has a finite solvable invariant series.
3) Decreasing chain of the commutators of the group are cut off on the single subgroup in finite number of steps.

***Statement 3***: The group $\mathbb{H}$ is solvable.

We prove that all the factors of the series $\mathbb{H} \geq A_i \geq A_e \geq \{e\}$ are Abelian.

Consider the factor group $\mathbb{H} / A_i$.

By Lagrange's theorem
$$|\mathbb{H}| = |\mathbb{H} : A_i| \cdot |A_i| \text{ then}$$
$$|\mathbb{H} : A_i| = 2 \text{ ie}$$
$\mathbb{H} / A_i$ is cyclic, and hence abelian.

Analogous, $|A_e : A_i| = 2$ it is $A_e / A_i$ is also an abelian group.

This means that
$$\mathbb{H} \geq A_i \geq A_e \geq \{e\} \text{ is a solvable series.}$$

Therefore $\mathbb{H}$ solvable group.

Since $[\mathbb{H}\mathbb{H}] = A_e$ then $\mathbb{H}$ the group is metabelian.

***Statement 4***: A group $G$ is metabelian if and only if its factor group by the center is also abelian.

We will prove this.

Let $h_1, h_2 \in \mathbb{H}$

$$h_1^{-1} h_2^{-1} h_1 h_2 \in A_e$$
$$h_1^{-1} h_2^{-1} \in h_2 h_1 A_e \quad \text{but}$$
$$h_1 h_2 \subset h_1 h_2 A_e \quad \text{then}$$
$$h_1 h_2 A_e = h_2 h_1 A_e$$

so $\mathbb{H} / A_i$ too is abelian.

The proof it is in the opposite direction we obtain the converse.

Since all the subgroups of group $\mathbb{H}$ are abelian groups then all the factor groups by the center are also abelian groups.

In the end we got:

The group $\mathbb{H}$ is not abelian, but all its subgroups are abelian groups.

All proper factor groups $\mathbb{H} / A_e$, $\mathbb{H} / A_i$ are abelian groups too.

All the factors $A_i / A_e$; $A_i / \{e\}$; $A_e / \{e\}$ are also abelian groups.

All proved properties are true for groups $A_j$, $A_k$.

***Definition***: Let $G$ a group under the operation of multiplication and is an arbitrary associative ring with unit 1.

Consider a set $K[G]$ consisting of all formal sums
$$\sum_{g \in G} \alpha_g g, \quad \alpha_g \in K$$
in which only a finite number of coefficients different from zero.

Formal sums $\sum_{g \in G} \alpha_g g$ and $\sum_{g \in G} \beta_g g$ ($\alpha_g, \beta_g \in K$) are equal if and only if.

If $x = \sum_{g \in G} \alpha_g g$ and $y = \sum_{g \in G} \beta_g g$ is elements of the set $K[G]$, then

$$x + y = \sum_{g \in G} (\alpha_g + \beta_g) g, \qquad xy = \sum_{g \in G} (\sum_{h \in G} \alpha_g \beta_{h^{-1}g}) g, \qquad \varepsilon = 1_K e_G$$

It is easy to see that a set $K[G]$ of respect to these operations is an associative ring with 1. This ring is called a group ring of the group $G$ over the ring $K$.

## Ideals of group ring $\mathbb{Z}_2[\mathbb{H}]$

Consider the ideals of ring $\mathbb{Z}_2[\mathbb{H}]$.
It is immediately verified that $\forall\ a, b \in \mathbb{H}$ next conditions are executed:
$$a^4 = b^4 = 1$$
$$a^2 = b^2 \qquad (1)$$
$$aba = b$$

Let $a := i,\ b := j$
Then

$$e = a^0 \qquad\qquad i = a$$
$$-e = a^2 \qquad\qquad -i = a^3$$

$$j = b \qquad\qquad k = ab$$
$$-j = a^2 b \qquad\qquad -k = a^3 b$$

The group ring $\mathbb{H}$ over the ring $\mathbb{Z}_2$, $\mathbb{Z}_2[\mathbb{H}]$ is consists of elements of the form
$$\alpha_{00} e + \alpha_{10} i + \alpha_{20}(-e) + \alpha_{30}(-i) + \alpha_{01} j + \alpha_{11} k + \alpha_{21}(-j) + \alpha_{31}(-k)$$
or if you use the notation
$$\alpha_{00} a^0 + \alpha_{10} a^1 + \alpha_{20} a^2 + \alpha_{30} a^3 + \alpha_{01} b + \alpha_{11} ab + \alpha_{21} a^2 b + \alpha_{31} a^3 b = \sum_{i=0}^{3} \sum_{j=0}^{i} \alpha_{ij} a^i b^j$$

This ring is not a field because it has zero divisors.
For example:
$$(e + a^2) = e + e = 0$$

We find the invertible elements.

***Statement 5***: Let $x = \sum_{i=0}^{3} \sum_{j=0}^{i} \alpha_{ij} a^i b^j$, $x \in \mathbb{H}$

Then, if

1) $\sum_{i=0}^{3} \sum_{j=0}^{i} \alpha_{ij} a^i b^j = 1$ in $\mathbb{Z}_2$

Then element $x$ is reversible.

2) $\sum_{i=0}^{3} \sum_{j=0}^{i} \alpha_{ij} a^i b^j = 0$ in $\mathbb{Z}_2$

Then element $x$ is a zero divisor.

From (1) and the fact $\alpha_{ij} \in \mathbb{Z}_2$ we get
$$(\sum_{i=0}^{3} \sum_{j=0}^{i} \alpha_{ij} a^i b^j)^4 = \sum_{i=0}^{3} \sum_{j=0}^{i} \alpha_{ij} (a^i)^4 (b^j)^4 = \sum_{i=0}^{3} \sum_{j=0}^{i} \alpha_{ij}$$

Hence we get that every element of $\mathbb{Z}_2[\mathbb{H}]$ in fourth degree is the sum of the coefficients $\alpha_{ij}$. That is, if the sequence $(\alpha_{00}, \alpha_{10}, ..., \alpha_{31})$ is an odd number of coefficients is different from $0_{\mathbb{Z}_2}$, then the corresponding element from $\mathbb{Z}_2[\mathbb{H}]$ in the fourth degree is $1_{\mathbb{Z}_2} e$ that is reversible.

And if even number of coefficients $\alpha_{ij} = 0_{\mathbb{Z}_2}$, then the element $x$ is a divisor of zero.

*Lemma 1*:
Let $G$ is a group,
$H$ is subgroup of group $G$,
$K$ is an arbitrary associative ring with 1.
The set of elements
$$L = \{u(h-1) \mid u \in R_l(G/H), 1 \neq h \in H\}$$
of group ring $K[G]$, is linearly independent over $K$ and forms a $K$-basis of module $I_l(H)$, where $R_l(G/H)$ is a complete system of representatives of left cosets of group $G$ by subgroup $H$ so that
$$G = \bigcup_{u \in R_l(G/h)} uH$$

$I_l(H)$ a set of elements of the group ring $K[G]$, which can be written as
$$\sum_{h \in H} x_h(h-1), \quad x_h \in K[G].$$

$I_l(H)$ a left ideal of the ring $K[G]$ and it is generated as a left ideal by elements
$$h-1, \quad h \in H.$$

*Proof:*
If the elements of the set $L$ are linearly dependent over $K$ then in the $L$ there is a minimal system of linearly independent elements
$$u_1(h_1 - 1), u_2(h_2 - 1), ......, u_s(h_s - 1)$$
Such that the relation
$$\lambda_1 u_1(h_1 - 1) + \lambda_2 u_2(h_2 - 1) + ... + \lambda_s u_s(h_s - 1) = 0 \qquad (2)$$
where $\lambda_1, \lambda_2, ... \lambda_s \in K$ and $\lambda_i \neq 0$, $u_i(h_i - 1) \in u_i H$.

Since different cosets of subgroup $H$ have no common elements, then (2) is possible when $u_1 = u_2 = ... = u_s$.
Then the elements $h_1, h_2, ..., h_s$ are different. From (2)
$$\lambda_1 h_1 + \lambda_2 h_2 + ... + \lambda_s h_s = \sum_{i=1}^{s} \lambda_i .$$
Which contradicts the linear independence $h_1, h_2, ..., h_s, 1$.
Consequently, the elements from $L$ are linearly independent.
Let
$$g \in G, h \in H \text{ и } g = uh_1 \ (u \in R_l(G/H); h_1(H)$$

then

$$g(h-1) = uh_1(h-1) = uh_1 h - uh_1 + u - u = u(h_1 h - 1) - u(h_1 - 1)$$
Hence the elements from $L$ form $K$-basis of $I_l(H)$.

*Proposal*:
If a subgroup $H$ is normal then $I_l(H)$ is a bilateral ideal.

*Proof:*

If $H$ is a normal in $G$, then $g^{-1}hg \in H \quad \forall g \in G, h \in H$
By virtue of
$$(h-1)g = g(g^{-1}hg - 1) \qquad g^{-1}hg - 1 \in H$$

$I_t(H)$ is a bilateral ideal.
Since all subgroups of group $\mathbb{H}$ are normal

$$\mathbb{H} \triangleright \begin{matrix} \{e,a,a^2,a^3\} \\ \{e,a^2,b,a^2b\} \\ \{e,a^2,ab,a^3b\} \end{matrix} \triangleright \{e,a^2\}$$

Then using Lemma 1, we get some ideals of the group ring $\mathbb{Z}_2[\mathbb{H}]$.

1) In the case a normal subgroup $H$ take the whole group $\mathbb{H}$. Then the basic elements for bilateral ideal $I(\mathbb{H})$ would be the following:

$$\begin{matrix}(e+a) & (e+b) & (e+a^3b) \\ (e+a^2) & (e+ab) & \\ (e+a^3) & (e+a^2b) & \end{matrix}$$

that is
$I(\mathbb{H}) = \{\lambda_1(e+a) + \lambda_2(e+a^2) + \lambda_3(e+a^3) + \lambda_4(e+b) + +\lambda_5(e+ab) + \lambda_6(e+a^2b) + \lambda_7(e+a^3b)\}$

$$\text{where } \lambda_1, \lambda_2, ..., \lambda_7 \in \mathbb{Z}_2$$

2) Consider the ideals of the relevant subgroups $<a>; <b>; <ab>$

$I(<a>) = \{\lambda_1(e+a) + \lambda_2(e+a^2) + \lambda_3(e+a^3) + \lambda_4(e+a)b + +\lambda_5(e+a^2)b + \lambda_7(e+a)^3 b\}$

That is, $H = <a>$, $R_r(\mathbb{H}/<a>) = \{e;b\}$, $\lambda_i \in Z_2$

Analogously we consider the ideals
$I(<b>)$,
$H = <b>$, $R_r(\mathbb{H}/<b>) = \{e;a\}$, $\lambda_i \in \mathbb{Z}_2$
$I(<ab>)$
$H = <ab>$, $R_r(\mathbb{H}/<ab>) = \{e;a\}$, $\lambda_i \in \mathbb{Z}_2$

3) $H = \{e;a^2\}$
Then $R_r\{\mathbb{H}/<a^2>\} = \{e;a;b;ab\}$,
$$I(<a^2>) = \{\lambda_1(e+a^2) + \lambda_2(e+a^2)b + \lambda_3(e+a^3)ab + \lambda_4(e+a^2)a\} \qquad \lambda_i \in \mathbb{Z}_2$$

Consider the other ideals of the group ring $\mathbb{Z}_2[\mathbb{H}]$.
The smallest ideal in the $\mathbb{Z}_2[\mathbb{H}]$ is the ideal which consist of two elements
$$\{e+a+a^2+a^3+b+ab+a^2b+a^3b;-0\}$$
It is immediately verified that the set
$$\{\lambda_1(e+a+a^2+a^3) + \lambda_2(b+ab+a^2b+a^3b)\} = I_a \qquad \lambda_i \in \mathbb{Z}_2$$
is an ideal, or
$$I_a = \{\lambda_1(e+a)^3 + \lambda_2(e+a)^3 b\} \qquad \lambda_i \in \mathbb{Z}_2$$
Analogously we consider the sets

$$I_b = \{\lambda_1 (e+b)^3 + \lambda_2 (e+b)^3 a\} \qquad \lambda_i \in \mathbb{Z}_2$$

and

$$I_{ab} = \{\lambda_1 (e+ab)^3 + \lambda_2 (e+ab)^3 a\} \qquad \lambda_i \in \mathbb{Z}_2$$

are ideals.

Combining bases of ideals
$I_a$ and $I_b$ $\{(e+a)^3; (e+a)^3 b; (e+b)^3; (e+b)^3 a\}$ and deleting linearly dependent elements we obtain a basis of the ideal which contains both $I_a$ and $I_b$.

$$(e+a)^3 + (e+a)^3 b = (e+b)^3 + (e+b)^3 a,$$

then the elements $(e+a)^3; (e+a)^3 b; (e+b)^3$ will be basis of ideal $I_{a \oplus b}$.

$$I_{a \oplus b} = \{\lambda_1 (e+a)^3 + \lambda_2 (e+a)^3 b + \lambda_3 (e+b)^3, \lambda_i \in \mathbb{Z}_2\}$$

In the same way we obtain the ideal $I_{a \oplus ab}; I_{b \oplus ab}$.

However $I_{a \oplus b} = I_{a \oplus ab} = I_{b \oplus ab}$.

It is immediately verified that the sets:

1) $\{(e+a)(e+a); (e+a)^2 a; (e+b)^2 b; (e+b)^2 ab\}$

2) $\{(e+a)(e+b); (e+a)(e+b)a; (e+a)(e+b)b; (e+a)(e+b)ab\}$

3) $\{(e+a)(e+ab); (e+a)(e+ab)a; (e+a)(e+ab)b; (e+a)(e+ab)ab\}$

are bases for the ideals that contain ideal $I_{a \oplus b}$.

We introduce the notation. Ideal corresponding to the set 1 denoted by a $M_a$

respectively 2 - $M_b$

3 - $M_b$.

Again, combining bases of ideals $M_a$ and $M_b$ and deleting linearly dependent elements we obtain a new basis

$$\{(e+a)^2; (e+a)^2 a; (e+b)^2 b; (e+b)^2 ab; (e+a)(e+b)\}$$

Ideal corresponding to this basis is denoted $M_{a \oplus b}$.

Analogously, combining bases of ideals $M_a$ and $M_{ab}$, $M_b$ and $M_{ab}$ and deleting linearly dependent elements, we obtain the ideals $M_{a \oplus ab}$ and $M_{b \oplus ab}$ respectively.

But

$$M_{a \oplus b} = M_{a \oplus ab} = M_{b \oplus ab}$$

Thus, using Lemma 1 and a direct check we got all the ideals $\mathbb{Z}_2[\mathbb{H}]$

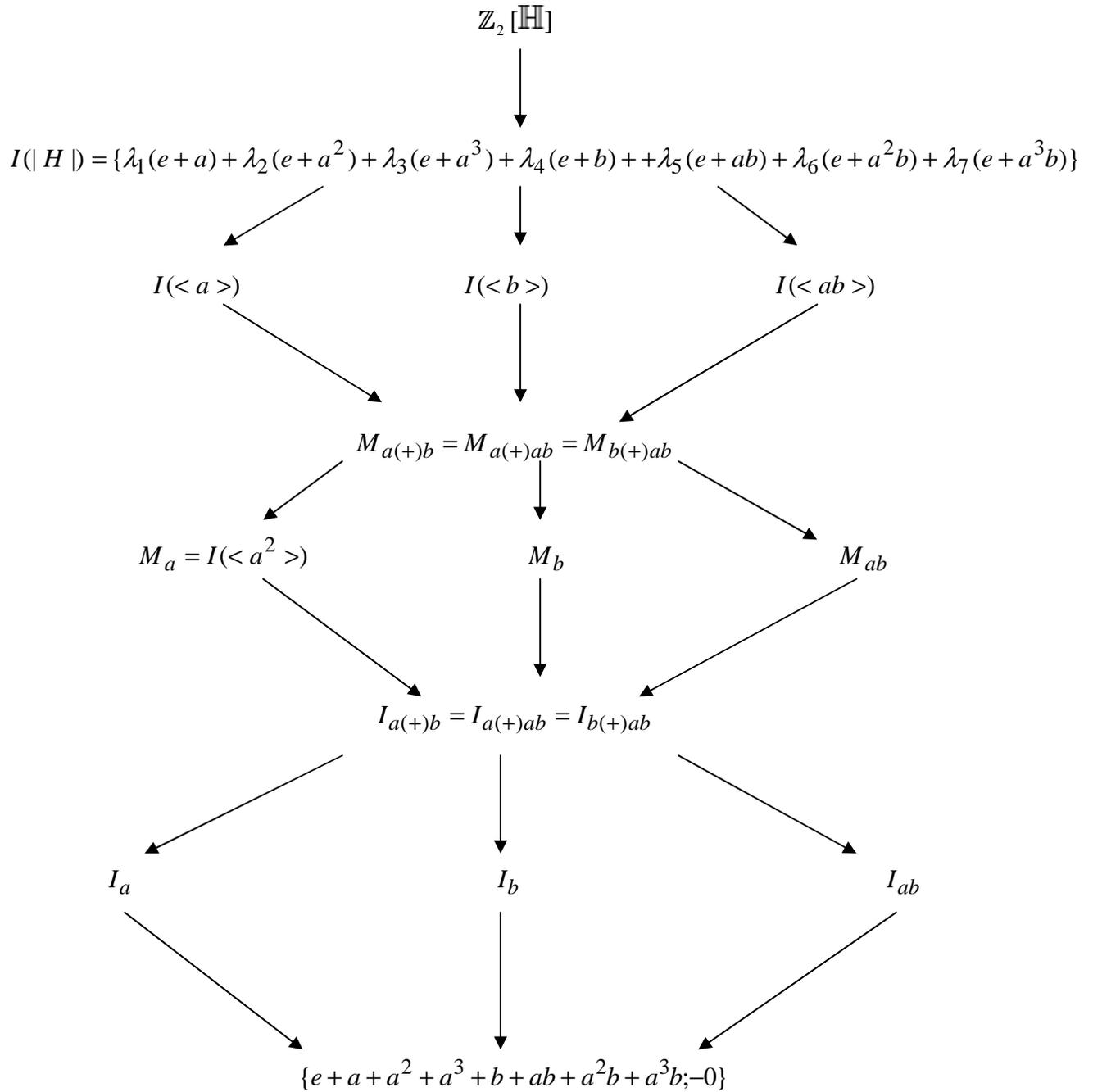

Consider the group ring $\mathbb{H}$ over the ring $\mathbb{Z}_3$, $\mathbb{Z}_3[\mathbb{H}]$ which are consists of elements of the form

$$\alpha_{00}e + \alpha_{10}i + \alpha_{20}(-e) + \alpha_{30}(-i) + \alpha_{01}j + \alpha_{11}k + \alpha_{21}(-j) + \alpha_{31}(-k)$$
$$\alpha_{ij} \in \mathbb{Z}_3, \ i = \overline{0,...,3}, \ j = 0,1$$

***Statement 6:***
Suppose $G$ is a finite group, $F$ is a field, $A$ is $FG$-module and $B$ is $FG$-submodule. If $\mathbf{\mathit{charF}} = 0$ or $\mathbf{\mathit{charF}} \notin \Pi(G)$ then there is a $FG$-submodule $C$ such that $A = B \oplus C$.
***Consequence***:
$A$ is $FG$-module, $G$ is a finite group, $F$ is a field.
If $\mathbf{\mathit{charF}} = 0$ or $\mathbf{\mathit{charF}} \notin \Pi(G)$ then $A = C_A(G) \oplus A(\omega FG)$

## Ideals of group ring $\mathbb{Z}_3[\mathbb{H}]$.

Consider the ideals of ring $\mathbb{Z}_3[\mathbb{H}]$.
The group ring $\mathbb{Z}_3[\mathbb{H}]$ satisfy the conditions of Statement 6. Then decompose $\mathbb{Z}_3[\mathbb{H}]$ into the direct sum of simple ideals, we get the ideals of the group ring $\mathbb{Z}_3[\mathbb{H}]$.
By Lemma 1, the ideal of a group ring is set generated by elements
$$h - e, \ h \in H, \ h \neq e$$
Thus, the basis an ideal are elements:

$$\begin{array}{lll} 2e + a & 2e + b & 2e + a^3b \\ 2e + a^2 & 2e + ab & \\ 2e + a^3 & 2e + a^2b & \end{array}$$

Ie
$$I(\mathbb{H}) = \{\lambda_1(2e + a) + \lambda_2(2e + a^2) + \lambda_3(2e + a^3) + \lambda_4(2e + b) +$$
$$+ \lambda_5(2e + ab) + \lambda_6(2e + a^2b) + \lambda_7(2e + a^3b)\}$$
$$\lambda_1, \lambda_2, ..., \lambda_7 \in \mathbb{Z}_3$$

By Statement 6 exists a submodule $C$ such that:
$$\mathbb{Z}_3[\mathbb{H}] = I(\mathbb{H}) \oplus C$$
In order to satisfy the condition $I(\mathbb{H}) \cap C = \varnothing$ necessary the bases of ideals $I(\mathbb{H})$ and $C$ are linearly independent.
As $C$ we can choose the set
$$\{\lambda(e + a + a^2 + a^3 + b + ab + a^2b + a^3b) : \lambda \in \mathbb{Z}_3\}$$
$C$ is an ideal in $\mathbb{Z}_3[\mathbb{H}]$.
Really
$$(\lambda(e + a + a^2 + a^3 + b + ab + a^2b + a^3b))^2 =$$
$$= 2\lambda^2(e + a + a^2 + a^3 + b + ab + a^2b + a^3b) \in C$$
Let
$$\lambda(e + a + a^2 + a^3 + b + ab + a^2b + a^3b)h = \lambda(e + a + a^2 + a^3 + b + ab + a^2b + a^3b)$$
Thus $C$ is an ideal.
We will prove that $I(\mathbb{H}) \cap C = \varnothing$.
Suppose the bases of $I(\mathbb{H})$ and $C$ are linearly dependent, then

$$\lambda(e+a+a^2+a^3+b+ab+a^2b+a^3b)h =$$
$$\alpha_1(2e+a)+\alpha_2(2e+a^2)+\alpha_3(2e+a^3)+\alpha_4(2e+b)+\alpha_5(2e+ab)+\alpha_6(2e+a^2b)+\alpha_7(2e+a^3b)$$

$$\begin{cases} \alpha_1 = \lambda \\ \alpha_2 = \lambda \\ \alpha_3 = \lambda \\ \alpha_4 = \lambda \\ \alpha_5 = \lambda \\ \alpha_6 = \lambda \\ \alpha_7 = \lambda \\ 2(\alpha_1 + \alpha_2 + ... + \alpha_7) = \lambda \end{cases}$$

such system has no solution in $\mathbb{Z}_3$. That is
$$I(\mathbb{H}) \cap C = \emptyset$$
Again $I(\mathbb{H})$ satisfies the conditions of Statement 6, then $I(\mathbb{H})$ can be decomposed into a direct sum of submodules.

The set generated by the basis

$$\begin{array}{ll} 2e+a & 2e+ab \\ 2e+a^2 & 2e+a^2b \\ 2e+a^3 & 2e+a^3b \end{array}$$

$$I(<a>) = \{\lambda_1(2e+a) + \lambda_2(2e+a^2) + \lambda_3(2e+a^3) + \lambda_4(2e+a)b + \lambda_5(2e+a^2)b + \lambda_7(2e+a)^3b\}$$
$$\lambda_1, \lambda_2, ..., \lambda_7 \in \mathbb{Z}_3$$

is an ideal of $I(\mathbb{H})$. That is $I(\mathbb{H}) = I(<a>) \oplus B$.

Let
$$B = \{\lambda(e+a+a^2+a^3+2b+2ab+2a^2b+2a^3b), \lambda \in \mathbb{Z}_3\}$$
It is immediately verified that $B$ is an ideal $I(\mathbb{H})$ and that $I(\mathbb{H}) \cap B = \emptyset$.

That is $\mathbb{Z}_3[\mathbb{H}]$ can be decomposed as follows:
$$\mathbb{Z}_3[\mathbb{H}] = C \oplus B \oplus I(<a>)$$

Again using Lemma 1 we obtain that the set
$$I(<a^2>) = \{\lambda_1(2e+a^2) + \lambda_2(2a+a^3) + \lambda_3(2b+a^2b) + \lambda_4(2ab+a^3b), \lambda_i \in \mathbb{Z}_3\}$$
is an ideal in $I(<a>)$.

that is
$$I(<a>) = I(<a^2>) \oplus D$$

in case $D$ we are selected set
$$D = \{\lambda_1(e+2a+a^2+2a^3) + \lambda_2(b+2ab+a^2b+2a^3b), \lambda_i \in \mathbb{Z}_3\}$$
Analogously we can verify that $D$ is an ideal in $I(<a>)$ and that $I(<a^2>) \cap D = \emptyset$.

In his turn $D$ decomposes into a direct sum of ideals
$$D = F \oplus G$$
where
$$F = \{\lambda(e+2a+a^2+2a^3+b+2ab+a^2b+2a^3b), \lambda \in \mathbb{Z}_3\}$$
$$G = \{\lambda(e+2a+a^2+2a^3+2b+ab+2a^2b+a^3b), \lambda \in \mathbb{Z}_3\}$$

Thus we have that:
$$\mathbb{Z}_3[\mathbb{H}] = C \oplus B \oplus F \oplus G \oplus I(<a^2>)$$

***Statement 7***:

The set $I(<a^2>)$ generated by the basis

$$2e + a^2 \qquad 2e + a^2 b$$
$$2e + a^3 \qquad 2e + a^3 b$$

contains no ideals.

To prove this, we use the Consequence of Statement 6.

The set $I(<a^2>)$ satisfies the Consequence. Then

$$I(<a^2>) = C_{I(<a^2>)}(<a^2>) \oplus I(<a^2>)(\omega \mathbb{Z}_3 <a^2>)$$

as

$$\mathbb{Z}_3 <a^2> = \{\lambda_1 e + \lambda_2 a^2, \ \lambda_i \in \mathbb{Z}_3\}$$

that

$$I(<a^2>)(\omega \mathbb{Z}_3 <a^2>) = I(<a^2>).$$

Thus we have

$$\mathbb{Z}_3[\mathbb{H}] =$$
$$\{\lambda_1(e + a + a^2 + a^3 + b + ab + a^2 b + a^3 b)\} \oplus$$

$$\oplus \{\lambda_2(e + a + a^2 + a^3 + 2b + 2ab + 2a^2 b + 2a^3 b)\} \oplus$$

$$\oplus \{\lambda_3(e + 2a + a^2 + 2a^3 + b + 2ab + a^2 b + 2a^3 b\} \oplus$$

$$\oplus \{\lambda_4(e + 2a + a^2 + 2a^3 + 2b + ab + 2a^2 b + a^3 b)\} \oplus$$

$$\oplus \{\alpha_1(2e + a^2) + \alpha_2(2a + a^3) + \alpha_3(2b + a^2 b) + \alpha_4(2ab + a^3 b)\}$$

$$\alpha_i, \lambda_j \in \mathbb{Z}_3.$$

Using this decomposition of $\mathbb{Z}_3[\mathbb{H}]$ we can get any ideal of $\mathbb{Z}_3[\mathbb{H}]$.